\documentclass[12pt]{article}%
\usepackage{graphicx}
\usepackage[intlimits]{amsmath}
\usepackage{latexsym}
\usepackage{amsfonts}
\usepackage{amssymb}%
 \makeatletter
\newfont{\footsc}{cmcsc10 at 8truept}
\newfont{\footbf}{cmbx10 at 8truept}
\newfont{\footrm}{cmr10 at 10truept}
\newtheorem{theorem}{Theorem}

\newtheorem{conjecture}[theorem]{Conjecture}

\newtheorem{lemma}[theorem]{Lemma}
\newtheorem{fact}[theorem]{Fact}

\newtheorem{proposition}[theorem]{Proposition}

\newenvironment{proof}[1][Proof]{\noindent{\textbf {#1}  }}  {\hfill$\Box$\bigskip}

\begin{document}

\title{The spectral radius of graphs without paths and cycles of specified length}
\author{Vladimir Nikiforov\\{\small Department of Mathematical Sciences, University of Memphis, Memphis TN
38152}}
\maketitle

\begin{abstract}
Let $G$ be a graph with $n$ vertices and $\mu\left(  G\right)  $ be the
largest eigenvalue of the adjacency matrix of $G.$ We study how large
$\mu\left(  G\right)  $ can be when $G$ does not contain cycles and paths of
specified order. In particular, we determine the maximum spectral radius of
graphs without paths of given length, and give tight bounds on the spectral
radius of graphs without given even cycles. We also raise a number of natural
open problems.\medskip


\end{abstract}

\section{Introduction}

How large can be the spectral radius $\mu\left(  G\right)  $ of a graph $G$ of
order $n$ without a path or cycle of order $k$? Such questions come easily
into mind when one studies the problems of extremal graph theory. Thus, recall
the general Tur\'{a}n type problem:\medskip

\emph{What is the maximum number of edges in a graph }$G$\emph{ of order }$n$
\emph{if }$G$\emph{ does not contain subgraphs of particular kind.}\medskip

In \cite{BrSo86}, Brualdi and Solheid raised an analogous spectral
problem:\medskip

\emph{What is }$\max$\emph{ }$\mu\left(  G\right)  $\emph{ if the graph }%
$G$\emph{ belongs to a specified class of graphs}.\medskip

Blending these two questions, we obtain a \emph{Brualdi-Solheid-Tur\'{a}n}
type problem:\medskip

\emph{What is }$\max$\emph{ }$\mu\left(  G\right)  $\emph{ if }$G$\emph{ is a
graph of order }$n$ \emph{and }$G$ \emph{does not contain subgraphs of
particular kind.}\medskip

Examples of such problems are numerous\ since to every Tur\'{a}n type problem
corresponds a Brualdi-Solheid-Tur\'{a}n type problem. In fact, many
fundamental types of graphs, like e.g. \emph{planar} or $k$\emph{-chromatic,}
are characterized by forbidden graphs, so the study of
Brualdi-Solheid-Tur\'{a}n type problems is an important topic in spectral
graph theory.

In this paper we focus on the maximum spectral radius of graphs of order $n$
without paths or cycles of specified length.

Write $C_{k}$ and $P_{k}$ for the cycle and path of order $k$, and let%
\begin{align*}
f_{l}\left(  n\right)   &  =\max\left\{  \mu\left(  G\right)  :\left\vert
G\right\vert =n,\text{ }C_{l}\nsubseteq G\right\}  ,\\
g_{l}\left(  n\right)   &  =\max\left\{  \mu\left(  G\right)  :\left\vert
G\right\vert =n,\text{ }C_{l}\nsubseteq G\text{ and }C_{l+1}\nsubseteq
G\right\}  \\
h_{l}\left(  n\right)   &  =\max\left\{  \mu\left(  G\right)  :\left\vert
G\right\vert =n,\text{ }P_{l}\nsubseteq G\right\}  .
\end{align*}
Perhaps it is more natural to define $g_{l}\left(  n\right)  $ as
$\max\left\{  \mu\left(  G\right)  :\left\vert G\right\vert =n,\text{ and
}C_{p}\nsubseteq G\text{ for }p\geq l\right\}  .$ However, this more relaxed
definition seems to determine precisely the same function $g_{l}\left(
n\right)  ,$ as suggested in Conjecture \ref{con1} in the concluding section
of this paper.

\subsubsection*{The value of $f_{l}\left(  n\right)  $ for odd $l$}

For odd $l$ the function $f_{l}\left(  n\right)  $ was essentially determined
in \cite{Nik08}: if $l$ is odd and $n>321l,$ then%
\[
f_{l}\left(  n\right)  =\sqrt{\left\lfloor n^{2}/4\right\rfloor }.
\]
The smallest ratio $n/l$ for which the equation is still valid is not known.
Note that the complete bipartite graph with color classes of size
$\left\lfloor n/2\right\rfloor $ and $\left\lceil n/2\right\rceil $ has no odd
cycles and its spectral radius is precisely $\sqrt{\left\lfloor n^{2}%
/4\right\rfloor }.$ Thus, for odd $l$ we have $f_{l}\left(  n\right)  \sim
n/2.$ As we shall see below, this is in sharp contrast with the value of
$f_{l}\left(  n\right)  $ for even $l.$

\subsubsection*{Bounds on $f_{4}\left(  n\right)  $}

The value of $f_{4}\left(  n\right)  $ was essentially determined in
\cite{Nik07}:\medskip

\emph{Let }$G$\emph{ be a graph of order }$n$\emph{ with }$\mu\left(
G\right)  =\mu$\emph{. If }$C_{4}\nsubseteq G,$\emph{ then}%
\begin{equation}
\mu^{2}-\mu\leq n-1. \label{c4in}%
\end{equation}
\emph{Equality holds if and only if every two vertices of }$G$\emph{ have
exactly one common neighbor, i.e., when }$G$\emph{ is the friendship
graph.}\medskip

An easy calculation implies that\emph{ }%
\[
f_{4}\left(  n\right)  =1/2+\sqrt{n-3/4}+O\left(  1/n\right)  ,
\]
where for odd $n$ the $O\left(  1/n\right)  $ term is zero. Finding the
precise value of $f_{4}\left(  n\right)  $ for even $n$ is an open problem.

\subsubsection*{Bounds on $f_{l}\left(  n\right)  $ for even $l>4$}

The inequality (\ref{c4in}) can be generalized for arbitrary even cycles in
the following way: if $C_{2k+2}\nsubseteq G,$ then%
\[
\mu^{2}-\left(  k-1\right)  \mu\leq k\left(  n-1\right)  .
\]
This inequality and a matching lower bound imply that%
\begin{equation}
\left(  k-1\right)  /2+\sqrt{kn}+o\left(  n\right)  \leq f_{2k+2}\left(
n\right)  \leq k/2+\sqrt{kn}+o\left(  n\right)  . \label{cyc}%
\end{equation}
The exact value of $f_{2k+2}\left(  n\right)  $ is not known for $k\geq2,$ and
finding this value seems a challenge. Nevertheless, the precision of
(\ref{cyc}) is somewhat surprising, given that the asymptotics of the maximum
number of edges in $C_{2k+2}$-free graphs of order $n$ is not known for
$k\geq2$.

\subsubsection*{Bounds on $g_{l}\left(  n\right)  $}

Let us consider now the function $g_{l}\left(  n\right)  $.\ To begin with,
Favaron, Mah\'{e}o, and Sacl\'{e} \cite{FMS93} showed that if a graph $G$ of
order $n$ contains neither $C_{3}$ nor $C_{4},$ then $\mu\left(  G\right)
\leq\sqrt{n-1}.$ Since the star of order $n$ has no cycles and its spectral
radius is $\sqrt{n-1},$ we deduce that
\[
g_{3}\left(  n\right)  =\sqrt{n-1}\emph{.}%
\]
We do not know the exact value of $g_{l}\left(  n\right)  $ for $l>3$.
Nevertheless, an example, together with the inequality $g_{2k+1}\left(
n\right)  \leq f_{2k+2}\left(  n\right)  $ and relation (\ref{cyc}), gives
\[
\left(  k-1\right)  /2+\sqrt{kn}+o\left(  n\right)  \leq g_{2k+1}\left(
n\right)  \leq k/2+\sqrt{kn}+o\left(  n\right)  ;
\]
thus, $g_{2k+1}\left(  n\right)  $ is known within an additive term not
exceeding $1/2.$

Luckily, for even $l$ we can give almost exact asymptotics of $g_{l}\left(
n\right)  :$
\[
g_{2k}\left(  n\right)  =\left(  k-1\right)  /2+\sqrt{kn}+O\left(
n^{-1/2}\right)  .
\]

\subsubsection*{Bounds on $h_{l}\left(  n\right)  $}

Finally, for $h_{l}\left(  n\right)  $ we have precise results when $n$ is
sufficiently large:
\begin{align*}
h_{2k}\left(  n\right)   &  =\left(  k-1\right)  /2+\sqrt{kn-\left(
3k^{2}+2k-1\right)  /4},\\
h_{2k+1}\left(  n\right)   &  =\left(  k-1\right)  /2+\sqrt{kn-\left(
3k^{2}+2k-1\right)  /4}+1/n+O\left(  n^{-3/2}\right)  .
\end{align*}
In addition, for every $l\geq4,$ we know the unique graph for which
$h_{l}\left(  n\right)  $ is attained when $n$ is sufficiently large.
Specifically, $h_{2k+1}\left(  n\right)  $ is known exactly but cannot be
given by a simple closed expression.

The main results of the paper are stated in the next section: first lower
bounds, and then upper bounds on $f_{2l}\left(  n\right)  ,$ $g_{l}\left(
n\right)  $ and $h_{l}\left(  n\right)  .$ The proofs of these results are
given in Section \ref{ps}. At the end of the paper we state two conjectures,
outlining possible solutions of related problems.

\section{Main results}

First we recall some notation, which in general follows \cite{Bol98}; thus, if
$G$ is a graph, we write:\medskip

- $V\left(  G\right)  $ for the vertex set of $G;$

- $\left\vert G\right\vert $ for the number of vertices of $G;$

- $E\left(  G\right)  $ for the edge set of $G$ and $e\left(  G\right)  $ for
$\left\vert E\left(  G\right)  \right\vert ;$

- $\delta\left(  G\right)  $ for the minimum degree of $G;$

- $G-u$ for the graph obtained by removing the vertex $u\in V\left(  G\right)
;$

- $\Gamma\left(  u\right)  $ for the set of neighbors of a vertex $u$ and
$d\left(  u\right)  $ for $\left\vert \Gamma\left(  u\right)  \right\vert ;$

- $e_{G}\left(  X\right)  $ for the number of edges induced by a set $X\subset
V\left(  G\right)  ;$

- $e_{G}\left(  X,Y\right)  $ for the number of edges joining vertices in $X$
to vertices in $Y,$ where $X$ and $Y$ are disjoint subsets of $V\left(
G\right)  $.\medskip

We write $K_{p}$ and $\overline{K}_{p}$ for the complete and the edgeless
graph of order $p$.\medskip

\subsection{Lower bounds on $f_{2l}\left(  n\right)  ,$ $g_{l}\left(
n\right)  $ and $h_{l}\left(  n\right)  $}

The lower bounds on $f_{2l}\left(  n\right)  ,$ $g_{l}\left(  n\right)  $ and
$h_{l}\left(  n\right)  $ are given by two families of graphs, which for
sufficiently large $n$ give the exact values of $h_{l}\left(  n\right)  $, and
perhaps also of $f_{2l}\left(  n\right)  $ and $g_{l}\left(  n\right)  ;$ for
a precise statement see Conjecture \ref{con1} in the concluding remarks.

Suppose that $1\leq k<n$ and let:

- $S_{n,k}$ be the graph obtained by joining every vertex of a complete graph
of order $k$ to every vertex of an independent set of order $n-k,$ that is to
say $S_{n,k}$ is the join of $K_{k}$ and $\overline{K}_{n-k};$

- $S_{n,k}^{+}$ be the graph obtained by adding one edge within the
independent set of $S_{n,k}.$

Clearly, $S_{n,k}$ and $S_{n,k}^{+}$ are graphs of order $n$ and%
\[
e\left(  S_{n,k}\right)  =kn-\left(  k^{2}+k\right)  /2,\text{ \ \ \ }e\left(
S_{n,k}^{+}\right)  =kn-\left(  k^{2}+k\right)  /2+1.
\]

To calculate $\mu\left(  S_{n,k}\right)  ,$ let $\mu=\mu\left(  S_{n,k}%
\right)  $ and apply a theorem of Finck and Grohmann \cite{FiGr65} (see also
\cite{CDS80}, Theorem 2.8) getting
\[
\mu^{2}-\left(  k-1\right)  \mu-k\left(  n-k\right)  =0.
\]
Thus, we have
\begin{equation}
\mu\left(  S_{n,k}\right)  =\left(  k-1\right)  /2+\sqrt{kn-\left(
3k^{2}+2k-1\right)  /4}. \label{musnk}%
\end{equation}

A slightly more involved approach gives $\mu\left(  S_{n,k}^{+}\right)  $ as follows.

\begin{proposition}
$\mu\left(  S_{n,k}^{+}\right)  $ is the largest root of the equation
\[
x^{3}-kx^{2}-\left(  kn-k^{2}-k+1\right)  x+k\left(  n-k-2\right)  =0,
\]
and satisfies the inequalities%
\begin{equation}
\frac{1}{n-k+\sqrt{kn/2}}<\mu\left(  S_{n,k}^{+}\right)  -\mu\left(
S_{n,k}\right)  <\frac{1}{n-k-2\sqrt{\left(  n-k\right)  /k}}. \label{bnds}%
\end{equation}

\end{proposition}

After some simple algebra, inequalities (\ref{musnk}) and (\ref{bnds}) give
\[
\mu\left(  S_{n,k}^{+}\right)  =\left(  k-1\right)  /2+\sqrt{kn-\left(
3k^{2}+2k-1\right)  /4}+1/n+O\left(  n^{-3/2}\right)  ,
\]

Note that $P_{l+1}\nsubseteq S_{n,k}$ and $C_{l}\nsubseteq S_{n,k}$ for
$l\geq2k+1.$ Likewise, $P_{l+1}\nsubseteq S_{n,k}$ and $C_{l}\nsubseteq
S_{n,k}$ for $l\geq2k+2.$ Therefore, we obtain the following bounds%
\begin{align}
h_{2k}\left(  n\right)   &  \geq\mu\left(  S_{n,k}\right)  =\left(
k-1\right)  /2+\sqrt{kn-\left(  3k^{2}+2k-1\right)  /4},\nonumber\\
h_{2k+1}\left(  n\right)   &  \geq\mu\left(  S_{n,k}^{+}\right)  =\left(
k-1\right)  /2+\sqrt{kn-\left(  3k^{2}+2k-1\right)  /4}+1/n+O\left(
n^{-3/2}\right)  ,\nonumber\\
g_{2k}\left(  n\right)   &  \geq\mu\left(  S_{n,k}\right)  =\left(
k-1\right)  /2+\sqrt{kn-\left(  3k^{2}+2k-1\right)  /4},\label{lbg}\\
g_{2k+1}\left(  n\right)   &  \geq\mu\left(  S_{n,k}^{+}\right)  =\left(
k-1\right)  /2+\sqrt{kn-\left(  3k^{2}+2k-1\right)  /4}+1/n+O\left(
n^{-3/2}\right)  ,\label{lbg1}\\
f_{2k+2}\left(  n\right)   &  \geq\mu\left(  S_{n,k}^{+}\right)  =\left(
k-1\right)  /2+\sqrt{kn-\left(  3k^{2}+2k-1\right)  /4}+1/n+O\left(
n^{-3/2}\right)  . \label{lbf}%
\end{align}

\subsection{Upper bounds on $f_{l}\left(  n\right)  ,$ $g_{l}\left(  n\right)
$ and $h_{l}\left(  n\right)  $}

\begin{theorem}
\label{th1} Let $k\geq1,$ $n\geq2^{4k}$ and $G$ be a graph of order $n.$

(a) If $\mu\left(  G\right)  \geq\mu\left(  S_{n,k}\right)  ,$ then $G$
contains a $P_{2k+2}$ unless $G=S_{n,k}.$

(b) If $\mu\left(  G\right)  \geq\mu\left(  S_{n,k}^{+}\right)  ,$ then $G$
contains a $P_{2k+3}$ unless $G=S_{n,k}^{+}.$
\end{theorem}

Theorem \ref{th1} implies that for every $k\geq1$ and $n\geq2^{4k},$ we have
\begin{align*}
h_{2k}\left(  n\right)   &  =\mu\left(  S_{n,k}\right)  =\left(  k-1\right)
/2+\sqrt{kn-\left(  3k^{2}+2k-1\right)  /4},\\
h_{2k+1}\left(  n\right)   &  =\mu\left(  S_{n,k}^{+}\right) \\
&  =\left(  k-1\right)  /2+\sqrt{kn-\left(  3k^{2}+2k-1\right)  /4}%
+1/n+O\left(  n^{-3/2}\right)  .
\end{align*}

\begin{theorem}
\label{th2} Let $k\geq1$ and $G$ be a graph of order $n.$ If
\begin{equation}
\mu\left(  G\right)  >k/2+\sqrt{kn+\left(  k^{2}-4k\right)  /4}, \label{lobo}%
\end{equation}
then $C_{2l+2}\subset G$ for every $l=1,\ldots,k$.
\end{theorem}

Theorem \ref{th2}, together with (\ref{lbg1}) and (\ref{lbf}), implies that%
\[
\left(  k-1\right)  /2+\sqrt{kn}+o\left(  n\right)  \leq g_{2k+1}\left(
n\right)  \leq f_{2k+2}\left(  n\right)  \leq k/2+\sqrt{kn}+o\left(  n\right)
.
\]
Finally, to determine the asymptotics of $g_{2k}\left(  n\right)  $ we need
the following theorem.

\begin{theorem}
\label{th3}Let $k\geq1$ and $G$ be a graph of order $n.$ If%
\[
\mu\left(  G\right)  >\left(  k-1\right)  /2+\sqrt{kn+\left(  k+1\right)
^{2}/4},
\]
then $C_{2k+1}\subset G$ or $C_{2k+2}\subset G.$
\end{theorem}

Theorem \ref{th3}, together with (\ref{lbg}), implies that for every $k\geq1,$%
\[
g_{2k}\left(  n\right)  =\left(  k-1\right)  /2+\sqrt{kn}+\Theta\left(
n^{-1/2}\right)  .
\]

\section{\label{ps}Proofs}

In this section we prove Theorems \ref{th1}, \ref{th2} and \ref{th3}.
Unfortunately, the proofs are involved and require a number of known facts and
some preliminary work.

\subsection{Some known facts}

We start with an upper bound on $\mu\left(  G\right)  $ that is particularly
efficient for our problems.

\begin{fact}
[\cite{Nik02}]\label{iNi}If $G$ is a graph with $n$ vertices, $m$ edges and
$\delta\left(  G\right)  =\delta,$ then
\begin{equation}
\mu\left(  G\right)  \leq\left(  \delta-1\right)  /2+\sqrt{2m-\delta n+\left(
\delta+1\right)  ^{2}/4}, \label{mainin}%
\end{equation}

\end{fact}

Note that for connected graphs inequality (\ref{mainin}) has been proved
independently by Hong, Shu and Fang \cite{HSF01}. A particular instance of
(\ref{mainin}) is the following upper bound (see \cite{Sta87}): \emph{for
every graph }$G$\emph{ with }$m$\emph{ edges,}%
\begin{equation}
\mu\left(  G\right)  \leq-1/2+\sqrt{2m+1/4}\leq\sqrt{2m}. \label{upb}%
\end{equation}

Next we turn to two classical results in extremal graph theory.

\begin{fact}
[Erd\H{o}s-Gallai (\cite{ErGa59}, Theorem 2.6)]\label{f1}Let $l\geq2$ and $G$
be a graph of order $n.$ If $e\left(  G\right)  >\left(  l/2\right)  n$, then
$G$ contains a $P_{l+2}.$
\end{fact}

Considerable work has been done to improve the above result under stronger
assumptions, see, e.g., \cite{BGLS08}, \cite{FaSc75} and \cite{Kop77}. In
particular, we shall use the following two facts, taken from \cite{BGLS08}.

\begin{fact}
\label{f2} Let $k\geq1,$ $n>3k$ and $G$ be a connected graph of order $n.$ If
\begin{equation}
e\left(  G\right)  \geq e\left(  S_{n,k}\right)  =kn-\left(  k^{2}+k\right)
/2, \label{maxed1}%
\end{equation}
then $G$ contains a $P_{2k+2},$ unless there is equality in (\ref{maxed1}) and
$G=S_{n,k}$
\end{fact}

\begin{fact}
\label{f3} Let $k\geq1,$ $n>3k$ and $G$ be a connected graph of order $n.$ If
\begin{equation}
e\left(  G\right)  \geq e\left(  S_{n,k}^{+}\right)  =kn-\left(
k^{2}+k\right)  /2+1, \label{maxed2}%
\end{equation}
then $G$ contains a $P_{2k+3},$ unless there is equality in (\ref{maxed2}) and
$G=S_{n,k}^{+}.$
\end{fact}

Recently in \cite{Nik08b}, Fact \ref{f1} has been extended as follows:

\begin{fact}
\label{lev6}Suppose that $k\geq1$ and let the vertices of a graph $G$ be
partitioned into two sets $U$ and $W$.

(A) If
\begin{equation}
2e_{G}\left(  U\right)  +e_{G}\left(  U,W\right)  >\left(  2k-2\right)
\left\vert U\right\vert +k\left\vert W\right\vert , \label{cond1}%
\end{equation}
then there exists a path of order $2k$ or $2k+1$ with both ends in $U.$

(B) If
\begin{equation}
2e_{G}\left(  U\right)  +e_{G}\left(  U,W\right)  >\left(  2k-1\right)
\left\vert U\right\vert +k\left\vert W\right\vert , \label{cond}%
\end{equation}
then there exists a path of order $2k+1$ with both ends in $U.$
\end{fact}

\subsection{Preliminary lemmas supporting the proof of Theorem \ref{th1}}

The proof of Theorems \ref{th1} is based on Lemmas \ref{lev1} through
\ref{lev4} below.

\begin{lemma}
\label{lev1}Let $G$ be a graph of order $n$ with minimum degree $\delta\left(
G\right)  =\delta$ and $\mu\left(  G\right)  =\mu.$ If $\left(  x_{1}%
,\ldots,x_{n}\right)  $ is a unit eigenvector to $\mu,$ then%
\[
\min\left\{  x_{1},\ldots,x_{n}\right\}  \leq\sqrt{\frac{\delta}{\mu
^{2}+\delta n-\delta^{2}}}.
\]

\end{lemma}

\begin{proof}
Set $\sigma=\min\left\{  x_{1},\ldots,x_{n}\right\}  .$ If $\sigma=0,$ the
assertion holds trivially, so let assume that $\sigma>0.$ This implies also
that $\delta>0$. Selecting $u\in V\left(  G\right)  $ to satisfy $d_{G}\left(
u\right)  =\delta,$ we have
\begin{align*}
\mu^{2}\sigma^{2}  &  \leq\mu^{2}x_{u}^{2}=\left(
{\textstyle\sum\limits_{i\in\Gamma\left(  u\right)  }}
x_{i}\right)  ^{2}\leq\delta%
{\textstyle\sum\limits_{i\in\Gamma\left(  u\right)  }}
x_{i}^{2}\leq\delta\left(  1-%
{\textstyle\sum\limits_{i\in V\left(  G\right)  \backslash\Gamma\left(
u\right)  }}
x_{i}^{2}\right) \\
&  \leq\delta\left(  1-\left(  n-\delta\right)  \sigma^{2}\right)
=\delta-\left(  \delta n-\delta^{2}\right)  \sigma^{2},
\end{align*}
implying that $\left(  \mu^{2}+\delta n-\delta^{2}\right)  \sigma^{2}%
\leq\delta.$ The desired inequality follows.
\end{proof}

\begin{lemma}
\label{lev2}Let $G$ be a graph of order $n$ and let $\left(  x_{1}%
,\ldots,x_{n}\right)  $ be a unit eigenvector to $\mu\left(  G\right)  .$ If
$u$ is\ a vertex satisfying $x_{u}=\min\left\{  x_{1},\ldots,x_{n}\right\}  ,$
then%
\[
\mu\left(  G-u\right)  \geq\mu\left(  G\right)  \frac{1-2x_{u}^{2}}%
{1-x_{u}^{2}}.
\]

\end{lemma}

\begin{proof}
Setting for short $\mu=\mu\left(  G\right)  ,$ we have%
\[
\mu x_{u}=\sum_{v\in\Gamma\left(  u\right)  }x_{v}\text{ \ \ and \ \ }%
\mu=2\sum_{vw\in E\left(  G\right)  }x_{v}x_{w}.
\]
Since Rayleigh's principle implies that%
\[
2\sum_{vw\in E\left(  G-u\right)  }x_{v}x_{w}\leq\mu\left(  G-u\right)
\sum_{v\in V\left(  G\right)  \backslash\left\{  u\right\}  }x_{v}^{2}%
=\mu\left(  G-u\right)  \left(  1-x_{u}^{2}\right)  ,
\]
we see that
\begin{align*}
\mu &  =2\sum_{vw\in E\left(  G-u\right)  }x_{v}x_{w}+2x_{u}\sum_{v\in
\Gamma\left(  u\right)  }x_{v}\\
&  =2\sum_{vw\in E\left(  G-u\right)  }x_{v}x_{w}+2x_{u}^{2}\mu\\
&  \leq\mu\left(  G-u\right)  \left(  1-x_{u}^{2}\right)  +2x_{u}^{2}\mu,
\end{align*}
and so,%
\[
\mu\left(  G-u\right)  \geq\mu\frac{1-2x_{u}^{2}}{1-x_{u}^{2}},
\]
as required.
\end{proof}

\begin{lemma}
\label{co1}Let $G$ be a graph of order $n,$ let $\mu\left(  G\right)  =\mu$
and $\left(  x_{1},\ldots,x_{n}\right)  $ be a unit eigenvector to $\mu.$ If
$u$ is\ a vertex satisfying $x_{u}=\min\left\{  x_{1},\ldots,x_{n}\right\}  ,$
then%
\[
\mu\left(  G-u\right)  \geq\mu\left(  1-\frac{1}{\mu^{2}/\delta+n-\delta
-1}\right)  .
\]

\end{lemma}

\begin{proof}
Lemma \ref{lev2} implies that%
\[
\mu\left(  G-u\right)  \geq\mu\frac{1-2x_{u}^{2}}{1-x_{u}^{2}}=\mu\left(
1-\frac{x_{u}^{2}}{1-x_{u}^{2}}\right)  .
\]
On the other hand, by Lemma \ref{lev1} we have%
\[
x_{u}^{2}\leq\frac{\delta}{\mu^{2}+\delta n-\delta^{2}},
\]
and so%
\begin{align*}
\mu\left(  G-u\right)   &  \geq\mu\left(  1-\frac{\delta}{\mu^{2}+\delta
n-\delta^{2}-\delta}\right) \\
&  =\mu\left(  1-\frac{1}{\mu^{2}/\delta+n-\delta-1}\right)  ,
\end{align*}
completing the proof.
\end{proof}

\begin{lemma}
\label{lev3}Let the numbers $a,$ $k,$ $n,$ $s$ satisfy%
\[
k\geq2,\text{ }s\geq1,\text{ and }n-s\geq4k^{3}+4\left\vert a\right\vert
\left(  k-1\right)  .
\]
Let the sequence $x_{0},\ldots,x_{s}$ satisfy%
\[
x_{0}\geq\left(  k-1\right)  /2+\sqrt{kn-a}%
\]
and%
\[
x_{i+1}\geq x_{i}\left(  1-\frac{1}{x_{i}^{2}/\left(  k-1\right)
+n-i-k}\right)  .
\]
for $0\leq i<s.$ Then, for every $i=1,\ldots,s,$ we have
\[
x_{i}\geq\left(  k-1\right)  /2+\sqrt{k\left(  n-i\right)  -a+1/2}.
\]

\end{lemma}

\begin{proof}
Clearly it is enough to prove the assertion for $i=1$ since it will follow by
induction for all $i=1,\ldots,s.$ Assume for a contradiction that
\[
x_{1}\leq\left(  k-1\right)  /2+\sqrt{k\left(  n-1\right)  -a+1/2}%
\]
and for short set
\[
b=\left(  k-1\right)  /2+\sqrt{kn-a}.
\]
Since the function
\[
x\left(  1-\frac{1}{x^{2}/\left(  k-1\right)  +n-k}\right)
\]
is increasing in $x,$ we have%
\begin{align*}
\left(  k-1\right)  /2+\sqrt{k\left(  n-1\right)  -a+1/2}  &  \geq x_{1}\\
&  \geq x_{0}\left(  1-\frac{1}{x_{0}^{2}/\left(  k-1\right)  +n-k}\right) \\
&  \geq b\left(  1-\frac{1}{\mu^{2}/\left(  k-1\right)  +n-k}\right)
\end{align*}
and so,
\begin{align*}
\frac{1}{b^{2}/\left(  k-1\right)  +n-k}  &  \geq\frac{b-\left(  \left(
k-1\right)  /2+\sqrt{k\left(  n-1\right)  -a+1/2}\right)  }{b}\\
&  =\frac{kn-a-\left(  k\left(  n-1\right)  -a+1/2\right)  }{b\left(
\sqrt{kn-a}+\sqrt{k\left(  n-1\right)  -a+1/2}\right)  }\\
&  =\frac{k-1/2}{b\left(  \sqrt{kn-a}+\sqrt{k\left(  n-1\right)
-a+1/2}\right)  }\\
&  >\frac{k-1/2}{2b\sqrt{kn-a}}.
\end{align*}
Hence, by the AM-GM inequality,%
\[
\frac{\sqrt{kn-a}}{k-1/2}>\frac{1}{2}\left(  \frac{b}{k-1}+\frac{n-k}%
{b}\right)  \geq\sqrt{\frac{n-k}{k-1}}.
\]
Squaring both sides of this inequality, we obtain%
\begin{align*}
n-k  &  <\left(  k-1\right)  \frac{kn-a}{\left(  k-1/2\right)  ^{2}}%
=\frac{\left(  k^{2}-k\right)  n}{k^{2}-k+1/4}-\frac{a\left(  k-1\right)
}{\left(  k-1/2\right)  ^{2}}\\
&  =n-\frac{n}{4\left(  k^{2}-k+1/4\right)  }-\frac{a\left(  k-1\right)
}{\left(  k-1/2\right)  ^{2}},
\end{align*}
and so,
\[
n<4\left(  k^{2}-k+1/4\right)  k-4a\left(  k-1\right)  <4k^{3}-4a\left(
k-1\right)  ,
\]
a contradiction completing the proof.
\end{proof}

\begin{lemma}
\label{lev4}Let the numbers $c\geq0,$ $k\geq2,$ $n\geq2^{4k},$ and let $G$ be
a graph of order $n.$ If $\delta\left(  G\right)  \leq k-1$ and%
\[
\mu\left(  G\right)  \geq\left(  k-1\right)  /2+\sqrt{kn-k^{2}+c},
\]
then there exists a graph $H$ satisfying one of the following conditions:

(i) $\mu\left(  H\right)  >\sqrt{\left(  2k+1\right)  \left\vert H\right\vert
};$

(ii) $\left\vert H\right\vert \geq\sqrt{n},$ $\delta\left(  H\right)  \geq k$
and
\[
\mu\left(  H\right)  >\left(  k-1\right)  /2+\sqrt{k\left\vert H\right\vert
-k^{2}+c+1/2}.
\]

\end{lemma}

\begin{proof}
Using the following procedure, define a sequence of graphs $G_{0},\ldots
,G_{k},$ satisfying $\left\vert G_{i}\right\vert =n-i$ for $i=0,\ldots,k:$

\qquad\textbf{begin}

\qquad\ \ \ \emph{set }$G_{0}=G;$

$\qquad\ \ \ \emph{set}$ $r=0;$

\qquad\ \ \ \textbf{while }$\mu\left(  H\right)  \leq\sqrt{\left(
2k+1\right)  \left\vert G_{r}\right\vert }$ and $\delta\left(  G_{r}\right)
\leq k-1$ \textbf{do }

\qquad\ \ \ \textbf{begin}

\qquad\ \ \ \ \ \ \ \emph{select a unit eigenvector }$\left(  x_{1}%
,\ldots,x_{n-r}\right)  $ \emph{to }$\mu\left(  G_{r}\right)  ;$

\qquad\ \ \ \ \ \ \ \emph{select a vertex }$u_{r}\in V\left(  G_{r}\right)  $
\emph{such that }$x_{u_{r}}=\min\left\{  x_{1},\ldots,x_{n-r}\right\}  ;$

\qquad\ \ \ \ \ \ \ \emph{set }$G_{r+1}=G_{r}-u_{r};$

\qquad\ \ \ \ \ \ \ \emph{add }$1$\emph{ to }$r;$

\qquad\ \ \ \textbf{end};

\qquad\textbf{end}.

Let $s=\min\left\{  r,n-\left\lfloor \sqrt{n}\right\rfloor \right\}  .$ Note
that for every $1\leq i<s,$ in view of $\delta\left(  G_{i}\right)  \leq k-1,$
Corollary \ref{co1} implies that%
\[
\mu\left(  G_{i+1}\right)  \geq\mu\left(  G_{i}\right)  \left(  1-\frac{1}%
{\mu^{2}\left(  G_{i}\right)  /\left(  k-1\right)  +n-i-k}\right)  .
\]
We shall prove that for every $i=1,\ldots,s,$
\begin{equation}
\mu\left(  G_{i}\right)  \geq\left(  k-1\right)  /2+\sqrt{k\left(  n-i\right)
-k^{2}+c+1/2} \label{in1}%
\end{equation}
Indeed, let $x_{i}=\mu\left(  G_{i}\right)  $ for $i=0,\ldots,s-1.$ Set
$a=k^{2}-c$ and note that
\[
\left\vert G_{i}\right\vert \geq\left\lfloor \sqrt{n}\right\rfloor \geq
2^{2k}\geq5k^{3}+\left\vert c\right\vert k\geq4k^{3}+\left\vert k^{2}%
-c\right\vert k\geq4k^{3}+\left\vert a\right\vert k.
\]
With this selection of $a,k,s,n,x_{0},\ldots,x_{s},$ Lemma \ref{lev3} implies
inequality (\ref{in1}).

Furthermore, for every $1\leq i<s,$ inequality (\ref{in1}) implies that%
\[
\mu^{2}\left(  G_{i}\right)  >k\left(  n-i-k\right)  ,
\]
and so, we find that%
\begin{align*}
\mu\left(  G_{i+1}\right)   &  \geq\mu\left(  G_{i}\right)  \left(  1-\frac
{1}{\mu^{2}\left(  G_{i}\right)  /\left(  k-1\right)  +n-i-k}\right) \\
&  >\mu\left(  G_{i}\right)  \left(  1-\frac{1}{k\left(  n-i-k\right)
/\left(  k-1\right)  +n-i-k}\right) \\
&  =\mu\left(  G_{i}\right)  \left(  1-\frac{k-1}{\left(  2k-1\right)  \left(
n-i-k\right)  }\right)  .
\end{align*}
On the other hand, Bernoulli's inequality gives%
\[
\left(  1-\alpha\frac{1}{i+1}\right)  \geq\frac{i^{\alpha}}{\left(
i+1\right)  ^{\alpha}}%
\]
whenever $0<\alpha<1$ and $i>0.$ In particular, we see that
\[
\left(  1-\frac{k-1}{\left(  2k-1\right)  \left(  n-i-k\right)  }\right)
\geq\left(  n-i+1-k\right)  ^{\left(  k-1\right)  /\left(  2k-1\right)
}\left(  n-i-k\right)  ^{-\left(  k-1\right)  /\left(  2k-1\right)  },
\]
and so, if $0\leq i<s,$ then%
\[
\mu\left(  G_{i+1}\right)  \left(  n-i+1-k\right)  ^{-\left(  k-1\right)
/\left(  2k-1\right)  }\geq\mu\left(  G_{i}\right)  \left(  n-i-k\right)
^{-\left(  k-1\right)  /\left(  2k-1\right)  }.
\]
Taking the first and last terms of this chain of inequalities, and setting
$H=G_{s},$ $p=\left\vert H\right\vert ,$ we find that
\begin{align*}
\mu\left(  H\right)  \left(  p-k\right)  ^{-\left(  k-1\right)  /\left(
2k-1\right)  }  &  \geq\mu\left(  G\right)  \left(  n-k\right)  ^{-\left(
k-1\right)  /\left(  2k-1\right)  }\\
&  >\sqrt{k}\left(  n-k\right)  ^{1/\left(  4k-2\right)  },
\end{align*}
and so,%
\begin{align*}
\mu\left(  H\right)   &  \geq\sqrt{k}\left(  n-k\right)  ^{1/\left(
4k-2\right)  }p^{\left(  k-1\right)  /\left(  2k-1\right)  }\\
&  =\sqrt{kp}\left(  n-k\right)  ^{1/\left(  4k-2\right)  }p^{-1/\left(
4k-2\right)  }.
\end{align*}

If $p=\left\lfloor \sqrt{n}\right\rfloor ,$ we see that $p\leq2\sqrt{n-k},$
and so,
\begin{align*}
\mu\left(  H\right)   &  \geq\sqrt{kp}\left(  n-k\right)  ^{1/\left(
4k-2\right)  }p^{-1/\left(  4k-2\right)  }\\
&  \geq\sqrt{kp}\left(  n-k\right)  ^{1/\left(  8k-4\right)  }2^{-1/\left(
4k-2\right)  }>\sqrt{\left(  2k+1\right)  p};
\end{align*}
thus $H$ satisfies condition \emph{(i)}$,$ completing the proof if
$p=\left\lfloor \sqrt{n}\right\rfloor $.

Otherwise, we have $p>\sqrt{n}$ and $\delta\left(  H\right)  \geq k.$ In view
of (\ref{in1}), $H$ satisfies condition \emph{(ii), }completing the proof of
Lemma \ref{lev4}.
\end{proof}

\subsubsection*{Proof of Theorem \ref{th1} for $k=1$}

For technical purposes we prove the case $k=1$ of Theorem \ref{th1}
separately. More precisely, we shall prove that if $G$ is a graph of order
$n,$ then:\medskip

\emph{(i) if }$n>5$\emph{ and }$\mu\left(  G\right)  \geq\sqrt{n-1},$\emph{
then }$G$\emph{ contains a }$P_{4},$\emph{ unless }$G=S_{n,1}.$

\emph{(ii) if }$n\geq10$\emph{ and }$\mu\left(  G\right)  \geq\mu\left(
S_{n,1}^{+}\right)  ,$\emph{ then }$G$\emph{ contains a }$P_{5},$\emph{ unless
}$G=S_{n,1}^{+}.$\medskip

\begin{proof}
[Proof of part (i)]First we first prove that every connected graph $G$ of
order at least $4$ satisfying $\mu\left(  G\right)  >\sqrt{n-1}$ contains a
$P_{4}$. Indeed, this is obvious if $G$ contains a $C_{4}.$ If $C_{4}%
\nsubseteq G,$ the result of Favaron, Mah\'{e}o, and Sacl\'{e} mentioned in
the Introduction implies that $G$ contains a triangle, say with vertices
$u,v,w.$ Since $G$ is connected, there is an edge between the sets $\left\{
u,v,w\right\}  $ and $V\left(  G\right)  \backslash\left\{  u,v,w\right\}  $,
so we get a $P_{4}.$

Suppose that $G$ is connected, $\mu\left(  G\right)  =\sqrt{n-1},$ and $G$
contains no $P_{4}.$ Then $G$ contains no cycles and so it must be a star
$S_{n,1}$.

Suppose now that $n>5$ and $G$ is not connected. Take a component $H$ with
$\mu\left(  H\right)  =\mu\left(  G\right)  \geq\sqrt{n-1}>\sqrt{\left\vert
H\right\vert -1}.$ According to the argument for connected graphs, $H$ must be
of order $3$ or $2.$ Thus, we have $2\geq\sqrt{n-1},$ a contradiction.\medskip

\begin{proof}
[Proof of part (ii)]First we shall prove that every connected graph $G$ of
order at least $5$ contains a $P_{5}$ unless $G$ has no cycles or
$G=S_{n,1}^{+}.$ Indeed, every cycle longer then $4$ contains $P_{5},$ so we
can assume that $G$ contains no such cycles.

If $G$ contains a $C_{4}$, say with vertices $u,v,w,t,$ there is an edge
between the sets $\left\{  u,v,w,t\right\}  $ and $V\left(  G\right)
\backslash\left\{  u,v,w,t\right\}  $, so we get a $P_{5}.$ Thus, we can
assume that $G$ contains no $C_{4}.$

Suppose $G$ contains a triangle, say with vertices $u,v,w.$ If two of the
vertices $\left\{  u,v,w\right\}  $ are joined to vertices from $V\left(
G\right)  \backslash\left\{  u,v,w\right\}  ,$ we get a $P_{5}.$ Hence only
one of the vertices $\left\{  u,v,w\right\}  $ is joined to vertices belonging
to $V\left(  G\right)  \backslash\left\{  u,v,w\right\}  ;$ let this be the
vertex $u$. Since all vertices $V\left(  G\right)  \backslash\left\{
u,v,w\right\}  $ are joined by some path to $u$, to avoid a $P_{5},$ all
vertices $V\left(  G\right)  \backslash\left\{  u,v,w\right\}  $ must be
joined to $u$ by an edge and the set $V\left(  G\right)  \backslash\left\{
u,v,w\right\}  $ must be independent. Therefore $G=S_{n,1}^{+}.$

Since a graph $G$ with no cycles satisfies $\mu\left(  G\right)  \leq
\sqrt{n-1}<\mu\left(  S_{n,1}^{+}\right)  ,$ part \emph{(ii)} is proved for
connected graphs of order at least $5$.

Let now $G$ be disconnected and $n\geq10$. Take a component $H$ with
$\mu\left(  H\right)  =\mu\left(  G\right)  \geq\mu\left(  S_{n,1}^{+}\right)
>\mu\left(  S_{\left\vert H\right\vert ,1}^{+}\right)  .$ According to the
argument for connected graphs, if $H$ contains no $P_{5},$ we have $\left\vert
H\right\vert \leq4.$ But this is impossible since $\mu\left(  H\right)
=\mu\left(  G\right)  >\sqrt{n-1}\geq3.$
\end{proof}
\end{proof}

\subsubsection*{Proof of Theorem \ref{th1} for $k\geq2$}

\begin{proof}
[ ]Clearly we can assume that $G$ is connected. For short set $m=e\left(
G\right)  .$\medskip

\emph{Proof of part (a).}

If $\delta\left(  G\right)  \geq k,$ inequality (\ref{mainin}) implies that
\begin{align*}
\mu\left(  G\right)   &  \leq\left(  \delta-1\right)  /2+\sqrt{2m-\delta
n+\left(  \delta+1\right)  ^{2}/4}\\
&  \leq\left(  k-1\right)  /2+\sqrt{2m-kn+\left(  k+1\right)  ^{2}/4}.
\end{align*}
Hence, in view of
\[
\mu\left(  G\right)  \geq\mu\left(  S_{n,k}\right)  =\left(  k-1\right)
/2+\sqrt{kn-k^{2}+\left(  k-1\right)  ^{2}/4},
\]
we obtain
\[
2m\geq2kn-\left(  k^{2}+k\right)  =2e\left(  S_{n,k}\right)  .
\]
Now Fact \ref{f2} implies that $G$ contains a $P_{2k+2}$ unless $G=P_{2k+2}$.
This completes the proof of part \emph{(a)} if $\delta\geq k.$

Assume now that $\delta\left(  G\right)  \leq k-1.$ Applying Lemma \ref{lev4}
with $c=\left(  k-1\right)  ^{2}/4$, we find a graph $H$ such that either
\emph{(i) }$\mu\left(  H\right)  >\sqrt{\left(  2k+1\right)  \left\vert
H\right\vert }$ or \emph{(ii)} $\left\vert H\right\vert >\sqrt{n},$
$\delta\left(  H\right)  \geq k$ and
\[
\mu\left(  H\right)  >\left(  k-1\right)  /2+\sqrt{k\left\vert H\right\vert
-k^{2}+\left(  k-1\right)  ^{2}/4}.
\]
If \emph{(i)} holds, then in view of (\ref{upb}), we see that
\[
2e\left(  H\right)  \geq\mu^{2}\left(  H\right)  >\left(  2k+1\right)
\left\vert H\right\vert ,
\]
and so, by Fact \ref{f1}, $G$ contains a $P_{2k+3},$ completing the proof of
part \emph{(a)} in this case\emph{.}

If \emph{(ii)} holds, then we have%
\[
\mu\left(  H\right)  >\left(  k-1\right)  /2+\sqrt{k\left\vert H\right\vert
-k^{2}+\left(  k-1\right)  ^{2}/4},
\]
and having from (\ref{mainin})%
\[
\mu\left(  H\right)  \leq\left(  k-1\right)  /2+\sqrt{2e\left(  H\right)
-k\left\vert H\right\vert +\left(  k-1\right)  ^{2}/4},
\]
we find that
\[
2e\left(  H\right)  >2kp-\left(  k^{2}+k\right)  =2e\left(  S_{n,k}\right)  .
\]
Fact \ref{f2} implies that $G$ contains a $P_{2k+2},$ completing the proof of
part \emph{(a)}.\medskip

\emph{Proof of part (b).}

The proof goes like in Part \emph{(a)}, but needs more care. If
\[
m\geq kn-\left(  k^{2}+k\right)  /2+1,
\]
the assertion follows from Fact \ref{f3}, so we shall assume that
\[
2m\leq2kn-\left(  k^{2}+k\right)  .
\]

If $\delta\left(  G\right)  \geq k,$ inequality (\ref{mainin}) gives%
\begin{align*}
\mu\left(  G\right)   &  \leq\left(  \delta-1\right)  /2+\sqrt{2m-\delta
n+\left(  \delta+1\right)  ^{2}/4}\\
&  \leq\left(  k-1\right)  /2+\sqrt{2m-kn+\left(  k+1\right)  ^{2}/4}\\
&  =\mu\left(  S_{n,k}\right)  <\mu\left(  S_{n,k}^{+}\right)  ,
\end{align*}
a contradiction. Thus we have $\delta\left(  G\right)  \leq k-1.$

Applying Lemma \ref{lev4} with $c=\left(  k-1\right)  ^{2}/4$, we find a graph
$H$ such that either \emph{(i) }$\mu\left(  H\right)  >\sqrt{\left(
2k+1\right)  \left\vert H\right\vert }$ or \emph{(ii)} $\left\vert
H\right\vert >\sqrt{n},$ $\delta\left(  H\right)  \geq k$ and
\[
\mu\left(  H\right)  >\left(  k-1\right)  /2+\sqrt{k\left\vert H\right\vert
-k^{2}+\left(  k-1\right)  ^{2}/4+1/2}.
\]
If \emph{(i)} holds, then in view of (\ref{upb}), we see that
\[
2e\left(  H\right)  \geq\mu^{2}\left(  H\right)  >\left(  2k+1\right)  p,
\]
and so, by Fact \ref{f1}, $G$ contains a $P_{2k+2}.$

If \emph{(ii)} holds, then
\begin{align*}
\mu\left(  H\right)   &  >\left(  k-1\right)  /2+\sqrt{k\left\vert
H\right\vert -k^{2}+\left(  k-1\right)  ^{2}/4+1/2}\\
&  =\left(  k-1\right)  /2+\sqrt{k\left\vert H\right\vert -\left(
3k^{2}+2k-3\right)  /4}.
\end{align*}
We shall show that
\[
\left(  k-1\right)  /2+\sqrt{k\left\vert H\right\vert -\left(  3k^{2}%
+2k-3\right)  /4}>\mu\left(  S_{p,k}^{+}\right)  .
\]
Indeed, assume for a contradiction that this inequality fails and set for
short $p=\left\vert H\right\vert $. In view of (\ref{bnds}), we see that
\begin{align*}
\mu\left(  S_{p,k}^{+}\right)   &  <\mu\left(  S_{p,k}\right)  +\frac
{1}{p-k-2\sqrt{\left(  p-k\right)  /k}}\\
&  =\left(  k-1\right)  /2+\sqrt{kp-\left(  3k^{2}+2k-1\right)  /4}+\frac
{1}{p-k-2\sqrt{\left(  p-k\right)  /k}}.
\end{align*}
Therefore,
\[
\sqrt{kp-\left(  3k^{2}+2k-3\right)  /4}-\sqrt{kp-\left(  3k^{2}+2k-1\right)
/4}<\frac{1}{p-k-2\sqrt{\left(  p-k\right)  /k}},
\]
and so,%
\begin{align*}
\frac{1}{p-k-2\sqrt{\left(  p-k\right)  /k}}  &  >\frac{1/2}{\sqrt{kp-\left(
3k^{2}+2k-3\right)  /4}+\sqrt{kp-\left(  3k^{2}+2k-1\right)  /4}}\\
&  >\frac{1}{4\sqrt{kp}}.
\end{align*}
Since $p\geq\left\lceil \sqrt{n}\right\rceil \geq2^{2k},$ the above inequality
is a contradiction.

The proof of Theorem \ref{th1} is completed.
\end{proof}

\subsubsection*{Proof of Theorem \ref{th2}}

\begin{proof}
[ ]Since the expression
\[
k/2+\sqrt{kn+\left(  k^{2}-4k\right)  /4}%
\]
is increasing in $k,$ it is enough to prove the existence only of $C_{2k+2}.$
Assume for a contradiction that $C_{2k+2}\nsubseteq G.$

Select $u\in V;$ let
\[
U=\Gamma_{G}\left(  u\right)  ,W=V\backslash\left(  \Gamma_{G}\left(
u\right)  \cup\left\{  u\right\}  \right)
\]
and set $H=G-u.$ Since $G$ contains no $C_{2k+2},$ $H$ contains no path of
order $2k+1$ whose ends belong to $U.$ By part (B) of Fact \ref{lev6}, we
have
\begin{align*}
2e_{H}\left(  U\right)  +e_{H}\left(  U,W\right)   &  \leq\left(  2k-1\right)
\left\vert A\right\vert +k\left\vert B\right\vert =\left(  2k-1\right)
d_{G}\left(  h\right)  +k\left(  n-d_{G}\left(  h\right)  -1\right) \\
&  =\left(  2k-1\right)  d_{G}\left(  h\right)  +k\left(  n-1\right)  .
\end{align*}
Hence, we see that%
\begin{align*}%
{\textstyle\sum\limits_{v\in U}}
d_{G}\left(  v\right)   &  =d_{G}\left(  u\right)  +%
{\textstyle\sum\limits_{v\in U}}
d_{H}\left(  v\right)  =d_{G}\left(  u\right)  +2e_{G}\left(  U\right)
+e_{G}\left(  U,W\right) \\
&  \leq d_{G}\left(  u\right)  +\left(  k-1\right)  d_{G}\left(  h\right)
+k\left(  n-1\right) \\
&  =kd_{G}\left(  u\right)  +k\left(  n-1\right)  .
\end{align*}
Letting $A$ be the adjacency matrix of $G,$ note that the $u$'th row sum of
the matrix
\[
A^{\prime}=A^{2}-kA
\]
is equal to
\[%
{\textstyle\sum\limits_{v\in\Gamma\left(  u\right)  }}
d_{G}\left(  v\right)  -kd_{G}\left(  u\right)  ,
\]
consequently, the maximum row sum $r_{\max}$ of $C$ satisfies%
\[
r_{\max}\leq k\left(  n-1\right)  .
\]
Letting $\mathbf{x}$ be an eigenvector of $A$ to $\mu,$ we see that the value
\[
\lambda=\mu^{2}-k\mu
\]
is an eigenvalue of $A^{\prime}$ with eigenvector $\mathbf{x}$. Therefore,
\[
\mu^{2}-k\mu=\lambda\leq k\left(  n-1\right)  ,
\]
and so,
\[
\mu\leq k/2+\sqrt{kn+\left(  k^{2}-4k\right)  /4}.
\]
This contradiction with (\ref{lobo}) completes the proof of Theorem \ref{th2}.
\end{proof}

\subsubsection*{Proof of Theorem \ref{th3}}

\begin{proof}
[ ]The proof of Theorem \ref{th3} is identical to the proof of Theorem
\ref{th2}, except that part (A) of Fact \ref{lev6} is used instead of part (B).
\end{proof}

\section{Two conjectures}

The following conjecture, if true, will give the exact values of the functions
$f_{2l}\left(  n\right)  $ and $g_{l}\left(  n\right)  $ for all $l>2$ and $n$
sufficiently large.

\begin{conjecture}
\label{con1}Let $k\geq2$ and $G$ be a graph of sufficiently large order $n.$

(a) if $\mu\left(  G\right)  \geq\mu\left(  S_{n,k}\right)  ,$ then $G$
contains $C_{2k+1}$ or $C_{2k+2}$ unless $G=S_{n,k};$

(b) if $\mu\left(  G\right)  \geq\mu\left(  S_{n,k}^{+}\right)  ,$ then $G$
contains $C_{2k+2}$ unless $G=S_{n,k}^{+}$.
\end{conjecture}

We finish with a conjecture that goes beyond cycles and paths. It is motivated
by the famous Erd\H{o}s-S\'{o}s conjecture about the maximum number of edges
in a graph of order $n$ that does not contain some tree of order $k.$

\begin{conjecture}
\label{con2}Let $k\geq2$ and $G$ be a graph of sufficiently large order $n.$

(a) if $\mu\left(  G\right)  \geq\mu\left(  S_{n,k}\right)  ,$ then $G$
contains all trees of order $2k+2$ unless $G=S_{n,k};$

(b) if $\mu\left(  G\right)  \geq\mu\left(  S_{n,k}^{+}\right)  ,$ then $G$
contains all trees of order $2k+3$ unless $G=S_{n,k}^{+}.$
\end{conjecture}

\end{document}